\documentclass[a4paper,12pt]{article}

\usepackage[english]{babel}
\usepackage{amsmath,amsfonts,amssymb,amsthm,bbm}

\usepackage[top=2cm, bottom=2cm, left=2.5cm, right=2.5cm]{geometry}

\usepackage{indentfirst}
\usepackage[T1]{fontenc} 

\usepackage{graphicx}
\usepackage{float}

\usepackage{epsfig}
\usepackage{epstopdf}

\numberwithin{equation}{section}

\renewcommand{\O}{\mathbb O}

\newcommand{\R}{\mathbb R}

\newcommand*{\abs}[1]{\left\lvert{#1}\right\rvert} 

\DeclareMathOperator*{\Ima}{\mathrm{Im}}
\DeclareMathOperator*{\Rea}{\mathrm{Re}}

\newcommand{\e}{\mathbf e}

\renewcommand{\i}{\mathbf i}
\renewcommand{\j}{\mathbf j}
\newcommand{\kk}{\mathbf k}

\usepackage{array}
\newcolumntype{L}[1]{>{\raggedright\let\newline\\\arraybackslash\hspace{0pt}}m{#1}}
\newcolumntype{C}[1]{>{\centering\let\newline\\\arraybackslash\hspace{0pt}}m{#1}}
\newcolumntype{R}[1]{>{\raggedleft\let\newline\\\arraybackslash\hspace{0pt}}m{#1}}

\begin{document}

\begin{center}
\Large
\textbf{Theoretical and numerical considerations \\ on the polar (exponential) form of octonions \\ and elements of higher-order Cayley-Dickson algebras}\\
\end{center}


\begin{center}
\begin{tabular}{cc}
	\textbf{{\L}ukasz B{\l}aszczyk\textsuperscript} \\
	l.blaszczyk@mini.pw.edu.pl 
\medskip \\
	Faculty of Mathematics and Information Science \\
	Warsaw University of Technology \\
	ul. Koszykowa 75, 00-662 Warszawa, Poland
\end{tabular}
\end{center}


\bigskip
\normalsize
\begin{abstract}
The article is devoted to the issue of the polar form of octonions. This is a~continuation of the works initiated by Hahn and Snopek in their articles from 2011. The results presented in the article show errors made in previous considerations and suggest the possibility of their improvement. The presented numerical method gives promising results, which as a result of further work can give analytical formulas for angles in the polar representation.
\end{abstract}

\section{Introduction}
The polar (exponential) representation is a convenient form of presenting complex numbers and has its own natural geometric interpretation. One of the applications that immediately follows from this representation is the use of complex numbers to the description of rotation on the $\R^2$ plane. In the case of hypercomplex numbers (especially quaternions and octonions) the situation is a bit more complicated.
\medskip

In our earlier works, we cited well-known results regarding polar representation (and consequently also trigonometric) of quaternions and octonions~\cite{BlaszczykPreludium2017}. The basic idea behind the representation of octonion in the polar form is the fact that any octonion \begin{align*}
o &= r_0 + r_1 \e_1 + r_2 \e_2 + r_3 \e_3 + r_4 \e_4 + r_5 \e_5 + r_6 \e_6 + r_7 \e_7 \in \O
\end{align*} can be rewritten as $o = \Rea o + \Ima o$, where $\Rea o = r_0$ is called the \emph{real part}, \linebreak and~$\Ima o = r_1 \e_1 + r_2 \e_2 + r_3 \e_3 + r_4 \e_4 + r_5 \e_5 + r_6 \e_6 + r_7 \e_7$ is the \emph{imaginary part}. Then, identically as in the case of complex numbers, we define the trigonometric form of every non-zero octonion $o\in\O$ as: \begin{align*}
o = \abs{o}\cdot(\cos\theta + \boldsymbol{\mu}\cdot\sin\theta),
\end{align*} where $\abs{o} = \sqrt{o\cdot o^*}$ is octonion norm, $\boldsymbol{\mu} = \frac{\Ima o}{\abs{\Ima o}}$ and $\theta\in\R$ is the solution of the system of equations \begin{align*}
\cos\theta = \frac{\Rea o}{\abs{o}},\qquad \sin\theta = \frac{\abs{\Ima o}}{\abs{o}}.
\end{align*} 
\medskip

Before we move from a trigonometric representation to an exponential form, we must first define an exponential function. Similarly as for the complex numbers and quaternions, we use the infinite series~\cite{Rodman_2014}. For any $o\in\O$, \begin{align*}
e^o := \sum\limits_{k=0}^{\infty} \frac{o^k}{k!}.
\end{align*} It can be shown that if we denote $\mathbf{o} = \Ima o$, then \begin{align*}
e^o = e^{\Rea o} \left(\cos \abs{\mathbf{o}} + \frac{\mathbf{o}}{\abs{\mathbf{o}}} \sin\abs{\mathbf{o}}\right).
\end{align*} We omit here a separate analysis of quaternions -- the reasoning is analogous and it suffices to notice that any quaternion is also an octonion, which has zero components standing at imaginary units $\e_4,\ldots,\e_7$. It should be noted that, due to the fact, that octonions are non-commutative, for any $o_1,o_2\in\O$ we have \begin{align*}
e^{o_1 + o_2} = e^{o_1}\cdot e^{o_2}\qquad\text{if and only if}\qquad o_1 \cdot o_2 = o_2\cdot o_1.
\end{align*} Using this notation we define the \emph{basic exponential form} of an octonion $o\in\O$, $o\neq 0$, as \begin{align} \label{eq:basic_exp}
o = \abs{o}\cdot e^{\theta\boldsymbol{\mu}},
\end{align} where $\theta$ and $\boldsymbol{\mu}$ are defined earlier. We can also write for any $\alpha\in\R$ that \begin{align*} 
\cos\alpha = \frac{1}{2}\left(e^{\boldsymbol{\mu}\alpha} + e^{-\boldsymbol{\mu}\alpha}\right), \qquad \sin\alpha = \frac{1}{2\boldsymbol{\mu}}\left(e^{\boldsymbol{\mu}\alpha} - e^{-\boldsymbol{\mu}\alpha}\right),
\end{align*} where $\boldsymbol{\mu}$ is any octonion such that $\abs{\boldsymbol{\mu}}=1$ and $\Rea\boldsymbol{\mu}=0$.
\medskip

While the exponential form given by the formula \eqref{eq:basic_exp} is convenient and shows that the hypercomplex numbers are a generalization of complex numbers, it does not give any more geometrical information. In the case of quaternions, this problem has already been solved. In his doctoral dissertation~\cite{bulow}, T. B\"{u}low showed that every non-zero quaternion can be presented in the form \begin{align*}
q = \abs{q} e^{\i\phi} e^{\kk\psi} e^{\j\theta},
\end{align*} where $(\phi,\psi,\theta)\in [-\pi,+\pi)\times [-\pi/2, +\pi/2) \times [-\pi/4,+\pi/4]$. The proof of the theorem on polar representation uses the algebraic properties of quaternion algebra and refers to the fact that each unitary and non-zero quaternion represents a certain rotation in $\R^3$ space. Angles $(2\phi,2\psi,2\theta)$ are known as Euler angles. In addition to stating the fact that such a representation is possible, B\"{u}low has given direct formulas that allow to calculate the values of these angles~\cite{bulow,Bulow2001}.
\medskip

In the case of higher order algebras (eg. octonions), the situation is much more complicated and so far no results have appeared in the literature showing complete evidence of polar form of octonions. In their works~\cite{HahnSnopek2016,HahnSnopek2011a,Snopek2013}, Hahn and Snopek presented a hypothesis regarding such representation, but the numerical tests they carried out showed that these are not correct formulas.
\medskip

In this article, we will present a different way of polar representation of octonions than the one presented in the works of Hahn and Snopek. In Section 2, we formulate the formulas presented in these works and comment on their erroneousness. Section 3 will be devoted to a relatively simple idea, which is based on the previously cited works, but avoids mistakes that were made. In Section 4, we will conduct numerical tests that will allow us to look at the obtained polar form of octonions and indicate the problems that are related to it. The report will end with a summary and discussion of the results in Section 5.
\medskip

\section{Past results presented in the literature}
The derivation of the formulas in the works of Hahn and Snopek was motivated by the analysis of hypercomplex analytical signals~\cite{HahnSnopek2016,HahnSnopek2011a,Snopek2013}. We can denote the octonion signal in the form \begin{align*}
o &= x_0 + x_1 \e_1 + x_2 \e_2 + x_3 \e_3 + x_4 \e_4 + x_5 \e_5 + x_6 \e_6 + x_7 \e_7
\end{align*} (in the original works it was a description of the analytical signal, whose OFT spectrum had the support only in the first octant of the $\R^3$ space). With this signal we can connect four complex signals (also analytical): \begin{align*}
u_0 &= (x_0 - x_3 - x_5 - x_6) + (x_1 + x_2 + x_4 - x_7)\i, \\
u_1 &= (x_0 + x_3 - x_5 + x_6) + (x_1 - x_2 + x_4 + x_7)\i, \\
u_2 &= (x_0 - x_3 + x_5 + x_6) + (x_1 + x_2 - x_4 + x_7)\i, \\
u_3 &= (x_0 + x_3 + x_5 - x_6) + (x_1 - x_2 - x_4 - x_7)\i,
\end{align*} which can be presented in exponential form: \begin{align*}
u_0 = \abs{u_0} e^{\i\varphi_0}, \qquad
u_1 = \abs{u_1} e^{\i\varphi_1}, \qquad
u_2 = \abs{u_2} e^{\i\varphi_2}, \qquad
u_3 = \abs{u_3} e^{\i\varphi_3}.
\end{align*} With the introduced notation, the authors proposed the following polar representation of octonion: \begin{align*}
o = \abs{o} e^{\e_1 \psi_1} \cdot e^{\e_3 \psi_3} \cdot e^{\e_2 \psi_2} \cdot e^{\e_7 \psi_7} \cdot e^{\e_4 \psi_4} \cdot e^{\e_6 \psi_6} \cdot e^{\e_5 \psi_5}.
\end{align*} The order of factors with successive imaginary units refers to the work of B\"{u}low.
\medskip

The individual angles values were proposed as follows: \begin{align*}
\psi_1 &= (\varphi_0 + \varphi_1 + \varphi_2 + \varphi_3)/4, \qquad
\psi_2 = (\varphi_0 + \varphi_1 - \varphi_2 - \varphi_3)/4, \\
\psi_4 &= (\varphi_0 - \varphi_1 + \varphi_2 - \varphi_3)/4, \qquad
\psi_5 = (\varphi_0 - \varphi_1 - \varphi_2 + \varphi_3)/4,
\end{align*} and the other angles are slightly more complicated: \begin{align*}
\sin(4\psi_3) &= \frac{\abs{u_0}^2 - \abs{u_1}^2}{\abs{u_0}^2 + \abs{u_1}^2}, \qquad
\sin(4\psi_6) = \frac{\abs{u_2}^2 - \abs{u_3}^2}{\abs{u_2}^2 + \abs{u_3}^2}, \\
\sin(4\psi_7) &= \frac{\abs{u_0}^2 + \abs{u_1}^2 - \abs{u_2}^2 - \abs{u_3}^3}{\abs{u_0}^2 + \abs{u_1}^2 + \abs{u_2}^2 + \abs{u_3}^2}.
\end{align*}
\medskip

\begin{figure}[!ht]
\centering \includegraphics[width = 0.8\textwidth]{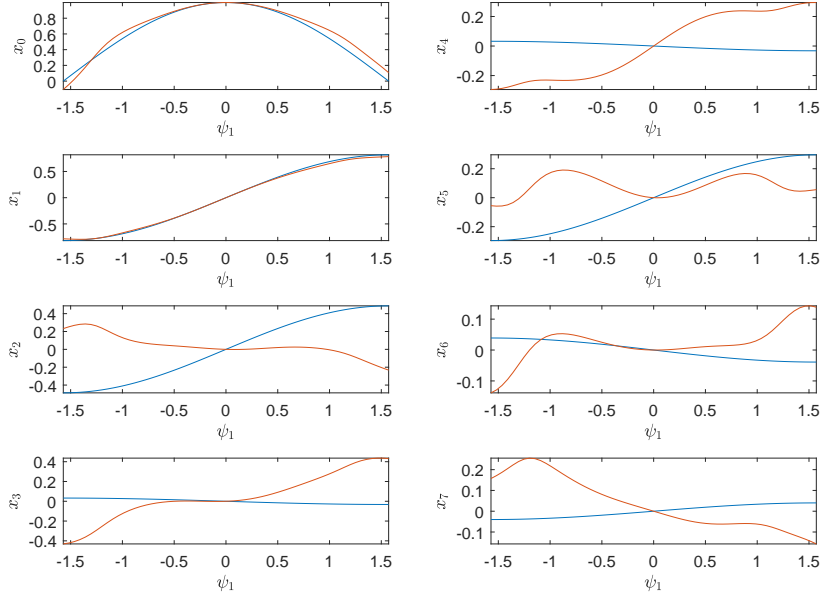}
\caption{Reconstructed coordinates $x_0,\ldots,x_7$ when changing parameter $\psi_1$ and fixed (random) parameters $\psi_2,\ldots,\psi_7$. The blue color corresponds to the correct value and the red value to the reconstructed value.}
\label{fig:a1}
\end{figure}

\begin{figure}[H]
\centering \includegraphics[width = 0.8\textwidth]{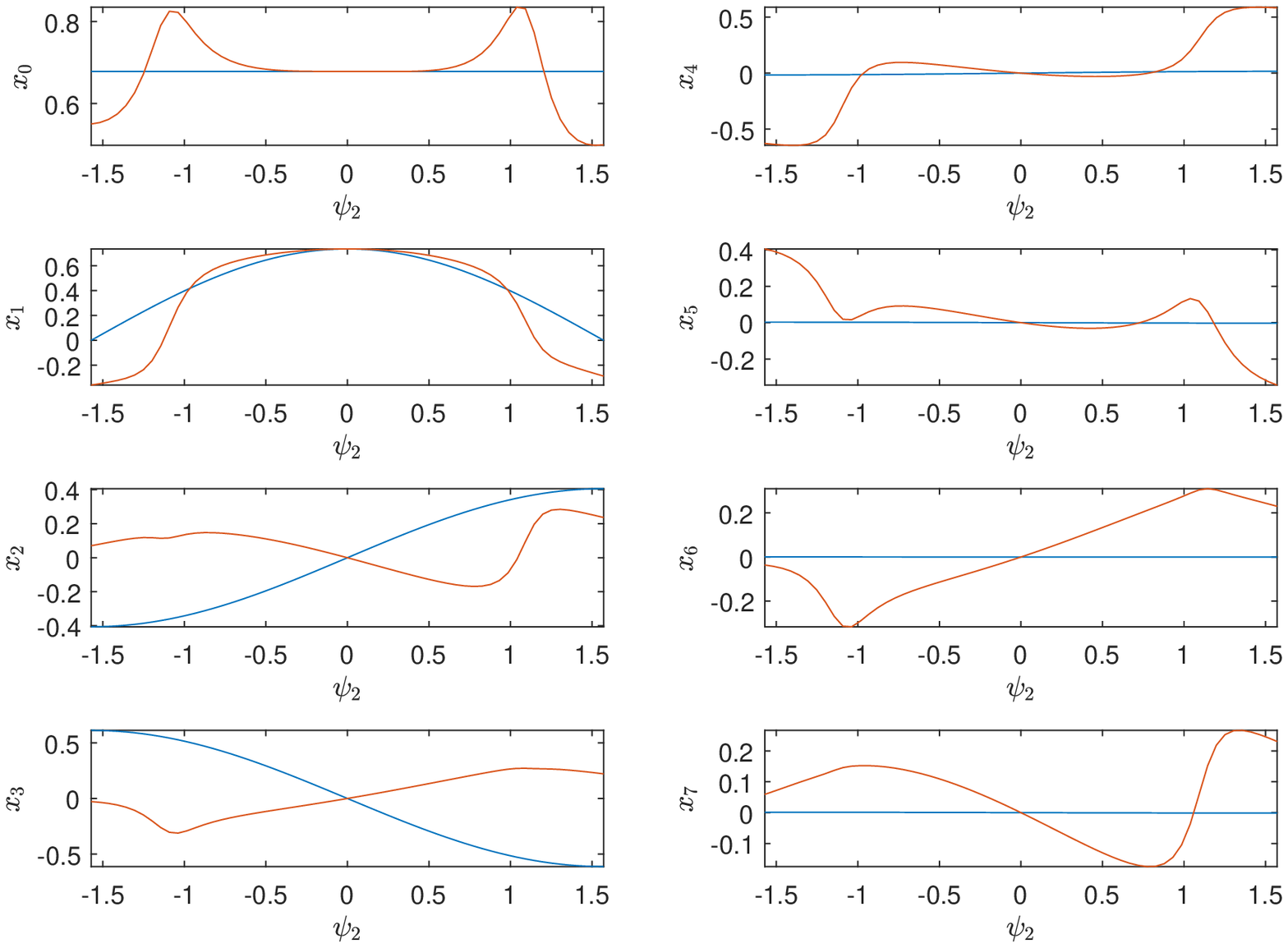}
\caption{Reconstructed coordinates $x_0,\ldots,x_7$ when changing parameter $\psi_2$ and fixed (random) parameters $\psi_1,\psi_3,\ldots,\psi_7$. The blue color corresponds to the correct value and the red value to the reconstructed value.}
\label{fig:a2}
\end{figure}

\begin{figure}[H]
\centering \includegraphics[width = 0.8\textwidth]{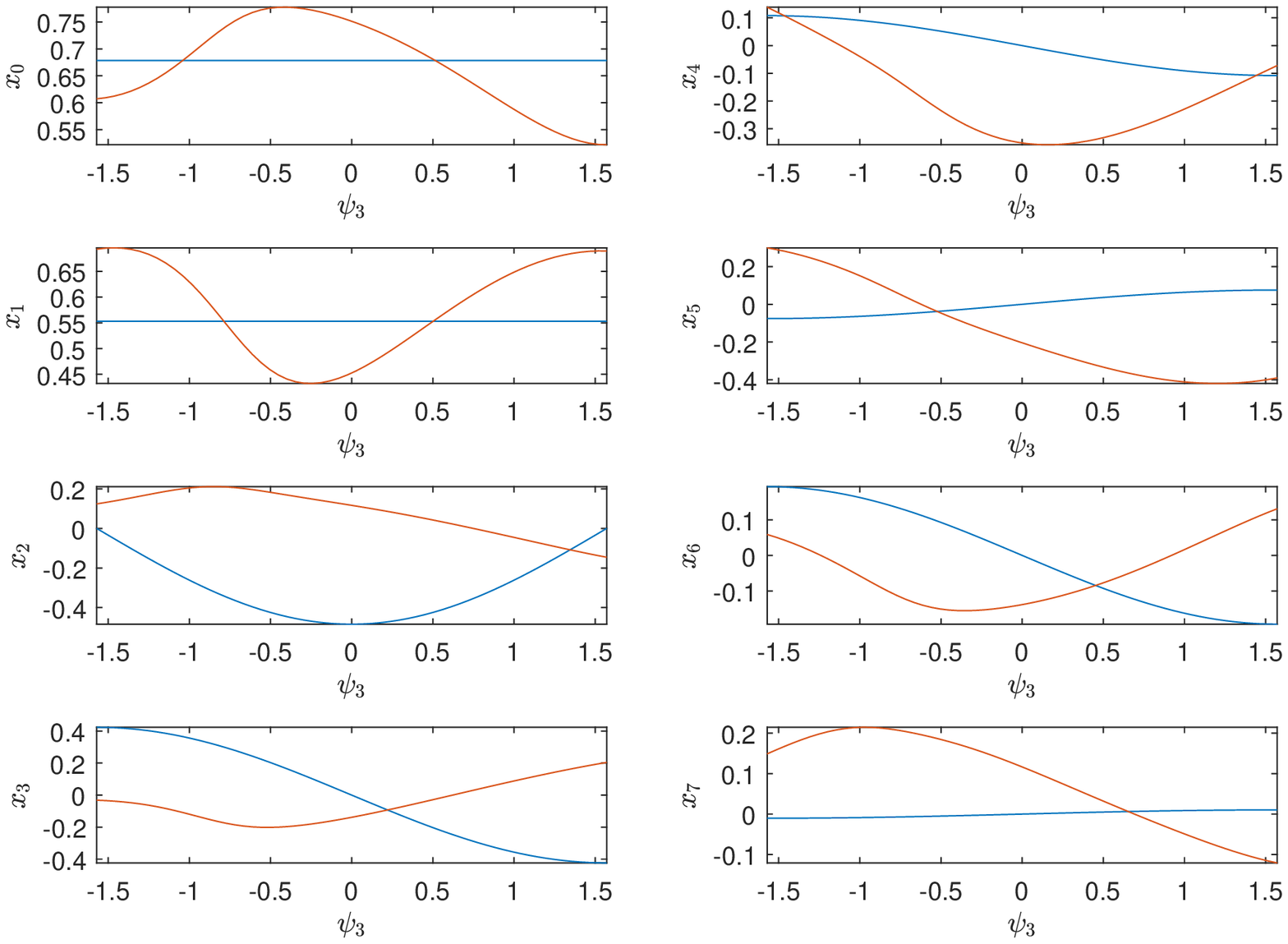}
\caption{Reconstructed coordinates $x_0,\ldots,x_7$ when changing parameter $\psi_3$ and fixed (random) parameters $\psi_1,\psi_2,\psi_4,\ldots,\psi_7$. The blue color corresponds to the correct value and the red value to the reconstructed value.}
\label{fig:a3}
\end{figure}

\begin{figure}[H]
\centering \includegraphics[width = 0.8\textwidth]{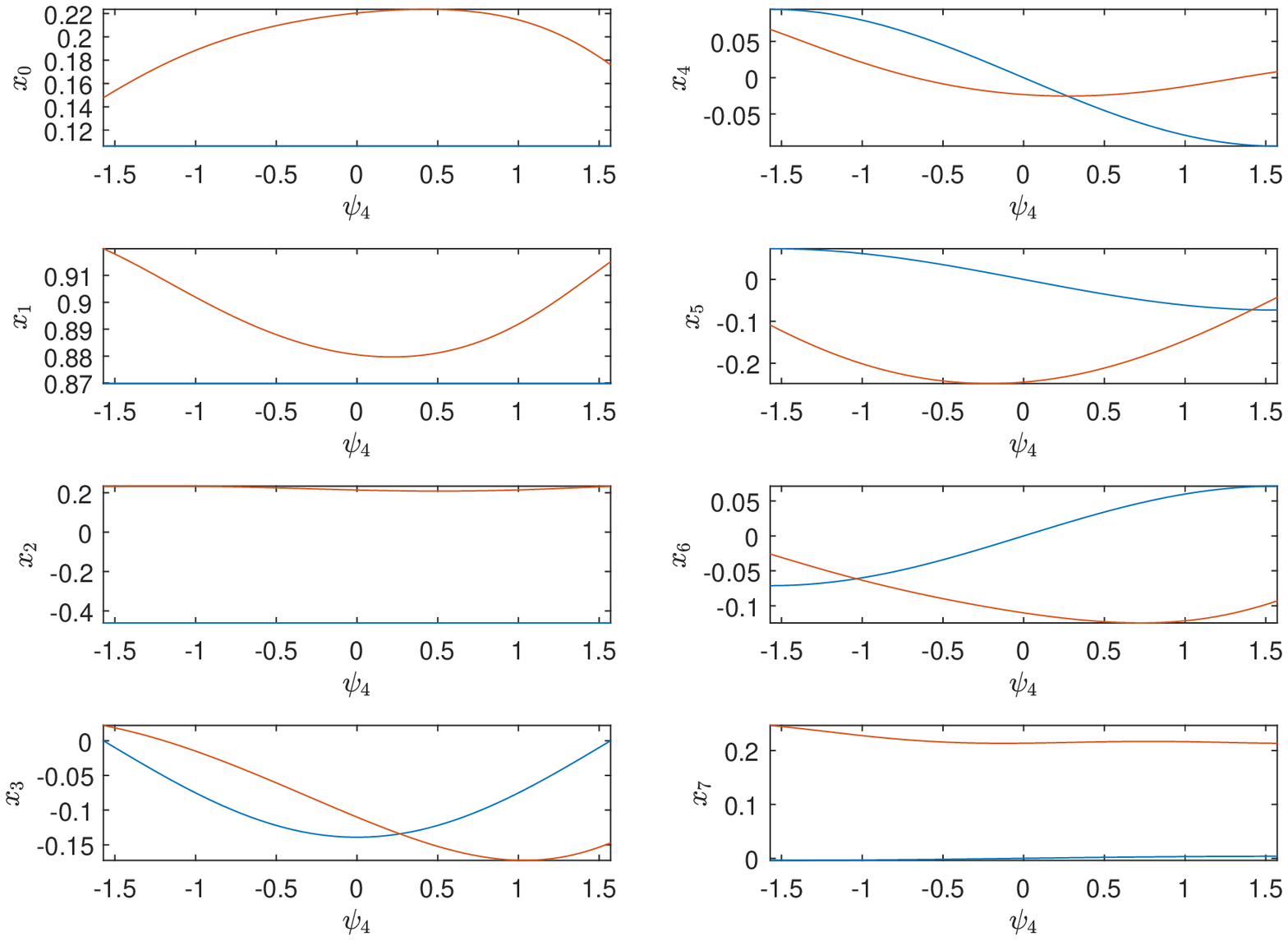}
\caption{Reconstructed coordinates $x_0,\ldots,x_7$ when changing parameter $\psi_4$ and fixed (random) parameters $\psi_1,\ldots,\psi_3,\psi_5,\ldots,\psi_7$. The blue color corresponds to the correct value and the red value to the reconstructed value.}
\label{fig:a4}
\end{figure}

\begin{figure}[H]
\centering \includegraphics[width = 0.8\textwidth]{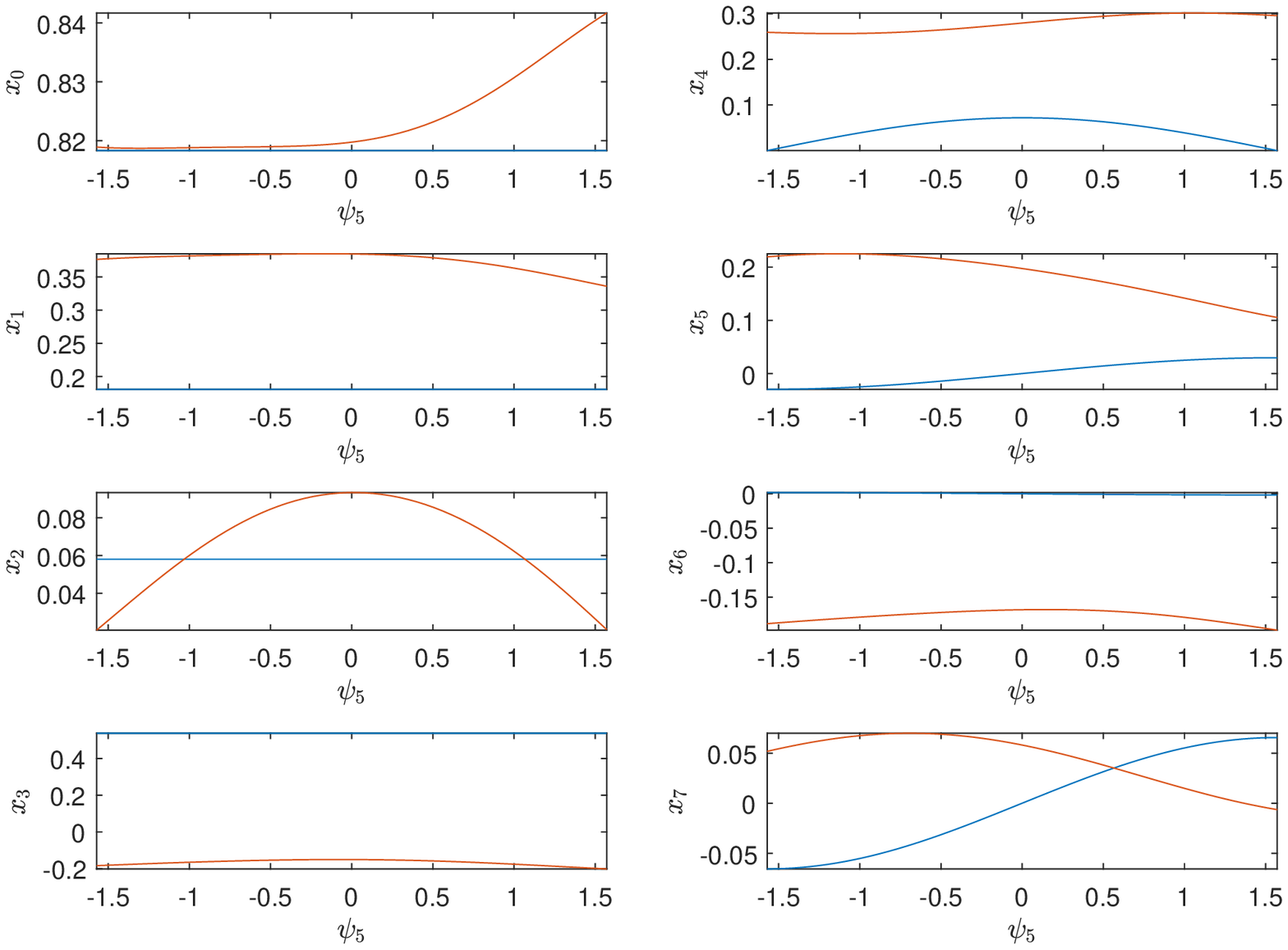}
\caption{Reconstructed coordinates $x_0,\ldots,x_7$ when changing parameter $\psi_5$ and fixed (random) parameters $\psi_1,\ldots,\psi_4,\psi_6,\psi_7$. The blue color corresponds to the correct value and the red value to the reconstructed value.}
\label{fig:a5}
\end{figure}

\begin{figure}[H]
\centering \includegraphics[width = 0.8\textwidth]{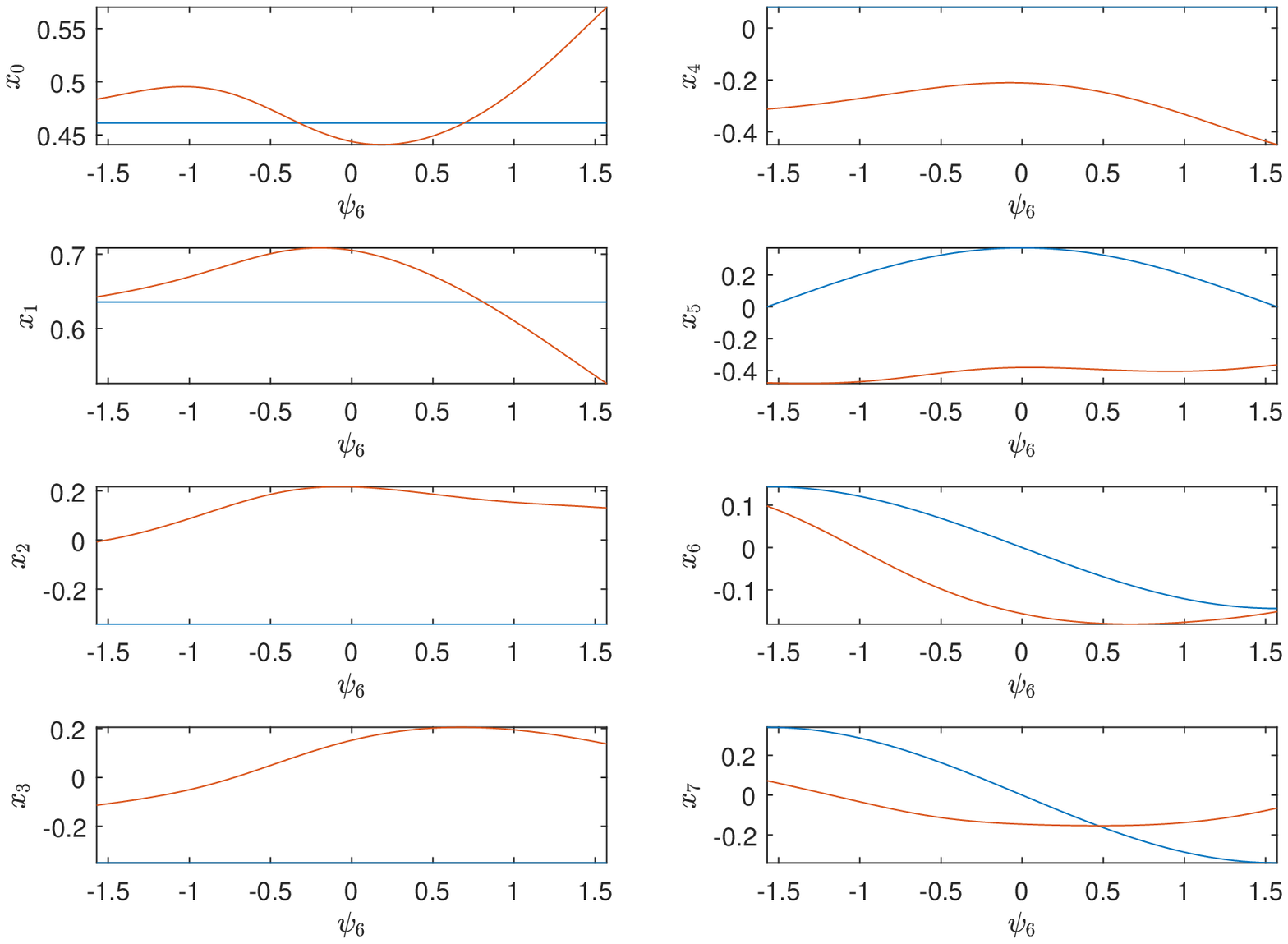}
\caption{Reconstructed coordinates $x_0,\ldots,x_7$ when changing parameter $\psi_6$ and fixed (random) parameters $\psi_1,\ldots,\psi_5,\psi_7$. The blue color corresponds to the correct value and the red value to the reconstructed value.}
\label{fig:a6}
\end{figure}

\begin{figure}[H]
\centering \includegraphics[width = 0.8\textwidth]{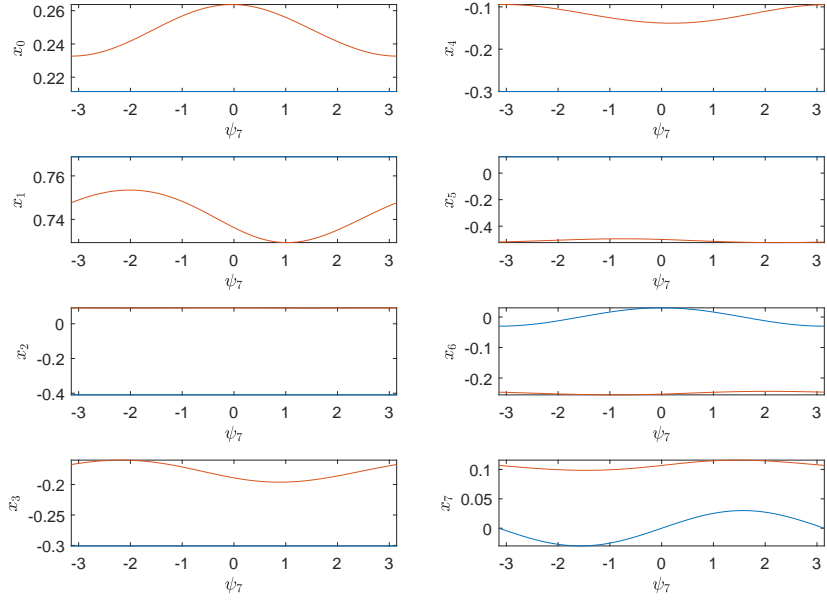}
\caption{Reconstructed coordinates $x_0,\ldots,x_7$ when changing parameter $\psi_7$ and fixed (random) parameters $\psi_1,\ldots,\psi_6$. The blue color corresponds to the correct value and the red value to the reconstructed value.}
\label{fig:a7}
\end{figure}

With the formulas given, it is easy to check their correctness. Since the set of unitary octonions can be identified with the unit sphere in $\R^8$, we generated a set of points on the basis of changing angles in spherical coordinates, i.e. numbers were generated as \begin{align*}
x_0 =& \cos\psi_1, \\
x_1 =& \sin\psi_1 \cos\psi_2, \\ 
x_2 =& \sin\psi_1 \sin\psi_2 \cos\psi_3, \\ 
x_3 =& \sin\psi_1 \sin\psi_2 \sin\psi_3 \cos\psi_4, \\ 
x_4 =& \sin\psi_1 \sin\psi_2 \sin\psi_3 \sin\psi_4 \cos\psi_5, \\ 
x_5 =& \sin\psi_1 \sin\psi_2 \sin\psi_3 \sin\psi_4 \sin\psi_5 \cos\psi_6, \\ 
x_6 =& \sin\psi_1 \sin\psi_2 \sin\psi_3 \sin\psi_4 \sin\psi_5 \sin\psi_6 \cos\psi_7, \\ 
x_7 =& \sin\psi_1 \sin\psi_2 \sin\psi_3 \sin\psi_4 \sin\psi_5 \sin\psi_6 \sin\psi_7, 
\end{align*} where $\psi_1,\ldots,\psi_6\in[-\pi/2,\pi/2)$ and $\psi_7\in[-\pi,\pi)$. During the experiment, we calculated the angles $\psi_1,\ldots,\psi_7$ using the previously given formulas and reconstructed all the coordinates of the octonion. In Fig. \ref{fig:a1}--\ref{fig:a7} the results of the experiment are shown. It can be easily seen that formulas are far from correct and only work in special cases.

\section{A simple idea}
We will focus on the factorization of any unit octonion $o\in \O$, i.e. \begin{align} \label{eq:forward_exp}
o = q\cdot e^{\e_4 \phi_4} \cdot e^{\e_5 \phi_5} \cdot e^{\e_6 \phi_6} \cdot e^{\e_7 \phi_7},
\end{align} where multiplication is performed from left to right, and $q$ is a certain unit quaternion. Note that at this stage, such a factorization is enough -- in his works B\"{u}low gave explicit formulas for the polar representation of a quaternion, i.e. $q = e^{\e_1 \phi_1} e^{\e_3 \phi_3} e^{\e_2 \phi_2}$~\cite{bulow,Bulow2001}.
\medskip

The polar form given by the formula \eqref{eq:forward_exp} is slightly different from that which was hipothesized in the works of Hahn and Snopek~\cite{HahnSnopek2011a}, but the way to derive it is so simple that it can be easily repeated to obtain a representation of the character \begin{align*}
o = q\cdot e^{\e_7 \psi_7} \cdot e^{\e_4 \psi_4} \cdot e^{\e_6 \psi_6} \cdot e^{\e_5 \psi_5}, \qquad
q = e^{\e_1 \psi_1} \cdot e^{\e_3 \psi_3} \cdot e^{\e_2 \psi_2}.
\end{align*}

Note that by doing a direct multiplication of particular factors and comparing real parts and individual imaginary parts, we get a system of equations \begin{align*}
x_0 =& \phantom{+}
         \cos(\phi_7)\cdot (\phantom{+}\cos(\phi_6)\cdot (\cos(\phi_5)\cdot \cos(\phi_4)\cdot y_0 - \sin(\phi_5)\cdot \sin(\phi_4)\cdot y_1) \\
&\phantom{+\cos(\phi_7)\cdot(} -\sin(\phi_6)\cdot (\cos(\phi_5)\cdot \sin(\phi_4)\cdot y_2 - \sin(\phi_5)\cdot \cos(\phi_4)\cdot y_3)) \\
&       -\sin(\phi_7)\cdot (-\sin(\phi_6)\cdot (\sin(\phi_5)\cdot \sin(\phi_4)\cdot y_0 + \cos(\phi_5)\cdot \cos(\phi_4)\cdot y_1) \\
&\phantom{+\cos(\phi_7)\cdot(} +\cos(\phi_6)\cdot (\sin(\phi_5)\cdot \cos(\phi_4)\cdot y_2 + \cos(\phi_5)\cdot \sin(\phi_4)\cdot y_3)), \displaybreak[3]\\
x_1 =& 
        -\sin(\phi_7)\cdot (\phantom{+}\sin(\phi_6)\cdot (\cos(\phi_5)\cdot \cos(\phi_4)\cdot y_0 - \sin(\phi_5)\cdot \sin(\phi_4)\cdot y_1) \\
&\phantom{+\cos(\phi_7)\cdot(} +\cos(\phi_6)\cdot (\cos(\phi_5)\cdot \sin(\phi_4)\cdot y_2 - \sin(\phi_5)\cdot \cos(\phi_4)\cdot y_3)) \\
&       +\cos(\phi_7)\cdot (\phantom{+}\cos(\phi_6)\cdot (\sin(\phi_5)\cdot \sin(\phi_4)\cdot y_0 + \cos(\phi_5)\cdot \cos(\phi_4)\cdot y_1) \\
&\phantom{+\cos(\phi_7)\cdot(} +\sin(\phi_6)\cdot (\sin(\phi_5)\cdot \cos(\phi_4)\cdot y_2 + \cos(\phi_5)\cdot \sin(\phi_4)\cdot y_3)), \displaybreak[3]\\
x_2 = & \phantom{+ }
         \sin(\phi_7)\cdot (\phantom{+}\cos(\phi_6)\cdot (\sin(\phi_5)\cdot \cos(\phi_4)\cdot y_0 + \cos(\phi_5)\cdot \sin(\phi_4)\cdot y_1) \\
&\phantom{+\cos(\phi_7)\cdot(} +\sin(\phi_6)\cdot (\sin(\phi_5)\cdot \sin(\phi_4)\cdot y_2 + \cos(\phi_5)\cdot \cos(\phi_4)\cdot y_3)) \\
&       +\cos(\phi_7)\cdot (\phantom{+}\sin(\phi_6)\cdot (\cos(\phi_5)\cdot \sin(\phi_4)\cdot y_0 - \sin(\phi_5)\cdot \cos(\phi_4)\cdot y_1) \\
&\phantom{+\cos(\phi_7)\cdot(} +\cos(\phi_6)\cdot (\cos(\phi_5)\cdot \cos(\phi_4)\cdot y_2 - \sin(\phi_5)\cdot \sin(\phi_4)\cdot y_3)), \displaybreak[3]\\
x_3 =& \phantom{+}
         \cos(\phi_7)\cdot (-\sin(\phi_6)\cdot (\sin(\phi_5)\cdot \cos(\phi_4)\cdot y_0 + \cos(\phi_5)\cdot \sin(\phi_4)\cdot y_1) \\
&\phantom{+\cos(\phi_7)\cdot(} +\cos(\phi_6)\cdot (\sin(\phi_5)\cdot \sin(\phi_4)\cdot y_2 + \cos(\phi_5)\cdot \cos(\phi_4)\cdot y_3)) \\
&       +\sin(\phi_7)\cdot (\phantom{+}\cos(\phi_6)\cdot (\cos(\phi_5)\cdot \sin(\phi_4)\cdot y_0 - \sin(\phi_5)\cdot \cos(\phi_4)\cdot y_1) \\
&\phantom{+\cos(\phi_7)\cdot(} -\sin(\phi_6)\cdot (\cos(\phi_5)\cdot \cos(\phi_4)\cdot y_2 - \sin(\phi_5)\cdot \sin(\phi_4)\cdot y_3)), \displaybreak[3]\\
x_4 =& 
        -\sin(\phi_7)\cdot (-\sin(\phi_6)\cdot (\sin(\phi_5)\cdot \cos(\phi_4)\cdot y_0 + \cos(\phi_5)\cdot \sin(\phi_4)\cdot y_1) \\
&\phantom{+\cos(\phi_7)\cdot(} +\cos(\phi_6)\cdot (\sin(\phi_5)\cdot \sin(\phi_4)\cdot y_2 + \cos(\phi_5)\cdot \cos(\phi_4)\cdot y_3)) \\
&       +\cos(\phi_7)\cdot (\phantom{+}\cos(\phi_6)\cdot (\cos(\phi_5)\cdot \sin(\phi_4)\cdot y_0 - \sin(\phi_5)\cdot \cos(\phi_4)\cdot y_1) \\
&\phantom{+\cos(\phi_7)\cdot(} -\sin(\phi_6)\cdot (\cos(\phi_5)\cdot \cos(\phi_4)\cdot y_2 - \sin(\phi_5)\cdot \sin(\phi_4)\cdot y_3)), \displaybreak[3]\\
x_5 =& \phantom{+}
         \cos(\phi_7)\cdot (\phantom{+}\cos(\phi_6)\cdot (\sin(\phi_5)\cdot \cos(\phi_4)\cdot y_0 + \cos(\phi_5)\cdot \sin(\phi_4)\cdot y_1) \\
&\phantom{+\cos(\phi_7)\cdot(} +\sin(\phi_6)\cdot (\sin(\phi_5)\cdot \sin(\phi_4)\cdot y_2 + \cos(\phi_5)\cdot \cos(\phi_4)\cdot y_3)) \\
&       -\sin(\phi_7)\cdot (\phantom{+}\sin(\phi_6)\cdot (\cos(\phi_5)\cdot \sin(\phi_4)\cdot y_0 - \sin(\phi_5)\cdot \cos(\phi_4)\cdot y_1) \\
&\phantom{+\cos(\phi_7)\cdot(} +\cos(\phi_6)\cdot (\cos(\phi_5)\cdot \cos(\phi_4)\cdot y_2 - \sin(\phi_5)\cdot \sin(\phi_4)\cdot y_3)), \displaybreak[3]\\
x_6 =& \phantom{+}
         \cos(\phi_7)\cdot (\phantom{+}\sin(\phi_6)\cdot (\cos(\phi_5)\cdot \cos(\phi_4)\cdot y_0 - \sin(\phi_5)\cdot \sin(\phi_4)\cdot y_1) \\
&\phantom{+\cos(\phi_7)\cdot(} +\cos(\phi_6)\cdot (\cos(\phi_5)\cdot \sin(\phi_4)\cdot y_2 - \sin(\phi_5)\cdot \cos(\phi_4)\cdot y_3)) \\
&       +\sin(\phi_7)\cdot (\phantom{+}\cos(\phi_6)\cdot (\sin(\phi_5)\cdot \sin(\phi_4)\cdot y_0 + \cos(\phi_5)\cdot \cos(\phi_4)\cdot y_1) \\
&\phantom{+\cos(\phi_7)\cdot(} +\sin(\phi_6)\cdot (\sin(\phi_5)\cdot \cos(\phi_4)\cdot y_2 + \cos(\phi_5)\cdot \sin(\phi_4)\cdot y_3)), \displaybreak[3]\\
x_7 =& \phantom{+}
         \sin(\phi_7)\cdot (\phantom{+}\cos(\phi_6)\cdot (\cos(\phi_5)\cdot \cos(\phi_4)\cdot y_0 - \sin(\phi_5)\cdot \sin(\phi_4)\cdot y_1) \\
&\phantom{+\cos(\phi_7)\cdot(} -\sin(\phi_6)\cdot (\cos(\phi_5)\cdot \sin(\phi_4)\cdot y_2 - \sin(\phi_5)\cdot \cos(\phi_4)\cdot y_3)) \\
&       +\cos(\phi_7)\cdot (-\sin(\phi_6)\cdot (\sin(\phi_5)\cdot \sin(\phi_4)\cdot y_0 + \cos(\phi_5)\cdot \cos(\phi_4)\cdot y_1) \\
&\phantom{+\cos(\phi_7)\cdot(} +\cos(\phi_6)\cdot (\sin(\phi_5)\cdot \cos(\phi_4)\cdot y_2 + \cos(\phi_5)\cdot \sin(\phi_4)\cdot y_3)),
\end{align*} where \begin{align} \label{eq:oct}
o &= x_0 + x_1 \e_1 + x_2 \e_2 + x_3 \e_3 + x_4 \e_4 + x_5 \e_5 + x_6 \e_6 + x_7 \e_7
\end{align} and \begin{align*}
q = y_0 + y_1 \e_1 + y_2 \e_2 + y_3 \e_3.
\end{align*} We obtained a system of eight non-linear algebraic equations with eight unknowns. At this point, one can try to solve this system of numerically, but one should not expect the uniqueness of such a solution. We will comment more on this fact in Section 4.
\medskip

As can be easily noticed, the expressions in individual equations are repeated. By introducing the appropriate notation, it can be shown that the determination of the polar form is reduced to carrying out a series of rotations in appropriate planes. Let's put: \begin{align*}
a_0 =& \phantom{+}\cos(\phi_6)\cdot (\cos(\phi_5)\cdot \cos(\phi_4)\cdot y_0 - \sin(\phi_5)\cdot \sin(\phi_4)\cdot y_1) \\
& -\sin(\phi_6)\cdot (\cos(\phi_5)\cdot \sin(\phi_4)\cdot y_2 - \sin(\phi_5)\cdot \cos(\phi_4)\cdot y_3) \displaybreak[3] \\
a_1 =& \phantom{+}\cos(\phi_6)\cdot (\sin(\phi_5)\cdot \sin(\phi_4)\cdot y_0 + \cos(\phi_5)\cdot \cos(\phi_4)\cdot y_1) \\
& +\sin(\phi_6)\cdot (\sin(\phi_5)\cdot \cos(\phi_4)\cdot y_2 + \cos(\phi_5)\cdot \sin(\phi_4)\cdot y_3) \displaybreak[3] \\
a_2 =& \phantom{+}\sin(\phi_6)\cdot (\cos(\phi_5)\cdot \sin(\phi_4)\cdot y_0 - \sin(\phi_5)\cdot \cos(\phi_4)\cdot y_1) \\
& +\cos(\phi_6)\cdot (\cos(\phi_5)\cdot \cos(\phi_4)\cdot y_2 - \sin(\phi_5)\cdot \sin(\phi_4)\cdot y_3) \displaybreak[3] \\
a_3 =& -\sin(\phi_6)\cdot (\sin(\phi_5)\cdot \cos(\phi_4)\cdot y_0 + \cos(\phi_5)\cdot \sin(\phi_4)\cdot y_1) \\
& +\cos(\phi_6)\cdot (\sin(\phi_5)\cdot \sin(\phi_4)\cdot y_2 + \cos(\phi_5)\cdot \cos(\phi_4)\cdot y_3) \displaybreak[3] \\
a_4 =& \phantom{+}\cos(\phi_6)\cdot (\cos(\phi_5)\cdot \sin(\phi_4)\cdot y_0 - \sin(\phi_5)\cdot \cos(\phi_4)\cdot y_1) \\
& -\sin(\phi_6)\cdot (\cos(\phi_5)\cdot \cos(\phi_4)\cdot y_2 - \sin(\phi_5)\cdot \sin(\phi_4)\cdot y_3) \displaybreak[3] \\
a_5 =& \phantom{+}\cos(\phi_6)\cdot (\sin(\phi_5)\cdot \cos(\phi_4)\cdot y_0 + \cos(\phi_5)\cdot \sin(\phi_4)\cdot y_1) \\
& +\sin(\phi_6)\cdot (\sin(\phi_5)\cdot \sin(\phi_4)\cdot y_2 + \cos(\phi_5)\cdot \cos(\phi_4)\cdot y_3) \displaybreak[3] \\
a_6 =& \phantom{+}\sin(\phi_6)\cdot (\cos(\phi_5)\cdot \cos(\phi_4)\cdot y_0 - \sin(\phi_5)\cdot \sin(\phi_4)\cdot y_1) \\
& +\cos(\phi_6)\cdot (\cos(\phi_5)\cdot \sin(\phi_4)\cdot y_2 - \sin(\phi_5)\cdot \cos(\phi_4)\cdot y_3) \displaybreak[3] \\
a_7 =& -\sin(\phi_6)\cdot (\sin(\phi_5)\cdot \sin(\phi_4)\cdot y_0 + \cos(\phi_5)\cdot \cos(\phi_4)\cdot y_1) \\
& +\cos(\phi_6)\cdot (\sin(\phi_5)\cdot \cos(\phi_4)\cdot y_2 + \cos(\phi_5)\cdot \sin(\phi_4)\cdot y_3).
\end{align*} Then \begin{align*}
& x_0 + x_7 \i = (a_0 + a_7 \i) e^{\i\phi_7}, \qquad
  x_1 + x_6 \i = (a_1 + a_6 \i) e^{\i\phi_7}, \\
& x_5 + x_2 \i = (a_5 + a_2 \i) e^{\i\phi_7}, \qquad
  x_4 + x_2 \i = (a_4 + a_2 \i) e^{\i\phi_7}.
\end{align*}
\medskip

Analogously, we follow the expressions in the definitions of numbers $a_0,\ldots,a_7$, obtaining \begin{align*}
b_0 =& \cos(\phi_5)\cdot \cos(\phi_4)\cdot y_0 - \sin(\phi_5)\cdot \sin(\phi_4)\cdot y_1 , \displaybreak[3] \\
b_1 =& \sin(\phi_5)\cdot \sin(\phi_4)\cdot y_0 + \cos(\phi_5)\cdot \cos(\phi_4)\cdot y_1 , \displaybreak[3] \\
b_2 =& \cos(\phi_5)\cdot \cos(\phi_4)\cdot y_2 - \sin(\phi_5)\cdot \sin(\phi_4)\cdot y_3 , \displaybreak[3] \\
b_3 =& \sin(\phi_5)\cdot \sin(\phi_4)\cdot y_2 + \cos(\phi_5)\cdot \cos(\phi_4)\cdot y_3 , \displaybreak[3] \\
b_4 =& \cos(\phi_5)\cdot \sin(\phi_4)\cdot y_0 - \sin(\phi_5)\cdot \cos(\phi_4)\cdot y_1 , \displaybreak[3] \\
b_5 =& \sin(\phi_5)\cdot \cos(\phi_4)\cdot y_0 + \cos(\phi_5)\cdot \sin(\phi_4)\cdot y_1 , \displaybreak[3] \\
b_6 =& \cos(\phi_5)\cdot \sin(\phi_4)\cdot y_2 - \sin(\phi_5)\cdot \cos(\phi_4)\cdot y_3 , \displaybreak[3] \\
b_7 =& \sin(\phi_5)\cdot \cos(\phi_4)\cdot y_2 + \cos(\phi_5)\cdot \sin(\phi_4)\cdot y_3 .
\end{align*} By re-conducting the previous reasoning we will get \begin{align*}
& a_0 + a_6 \i = (b_0 + b_6 \i) e^{\i\phi_6}, \qquad
  a_7 + a_1 \i = (b_7 + b_1 \i) e^{\i\phi_6}, \\
& a_4 + a_2 \i = (b_4 + b_2 \i) e^{\i\phi_6}, \qquad
  a_3 + a_5 \i = (b_3 + b_5 \i) e^{\i\phi_6}.
\end{align*} The next step involves the quaternion coefficients $q$, i.e. \begin{align*}
c_0 =& \cos(\phi_4)\cdot y_0, \qquad c_4 = \sin(\phi_4)\cdot y_0, \qquad
c_1 = \cos(\phi_4)\cdot y_1, \qquad c_5 = \sin(\phi_4)\cdot y_1, \\
c_2 =& \cos(\phi_4)\cdot y_2, \qquad c_6 = \sin(\phi_4)\cdot y_2, \qquad
c_3 = \cos(\phi_4)\cdot y_3, \qquad c_7 = \sin(\phi_4)\cdot y_3, 
\end{align*} and then \begin{align*}
& b_0 + b_5 \i = (c_0 + c_5 \i) e^{\i\phi_5}, \qquad
  b_4 + b_1 \i = (c_4 + c_1 \i) e^{\i\phi_5}, \\
& b_2 + b_7 \i = (c_2 + c_7 \i) e^{\i\phi_5}, \qquad
  b_6 + b_3 \i = (c_6 + c_3 \i) e^{\i\phi_5}.
\end{align*} Ultimately, using the above notation, we will get \begin{align*}
c_0 + c_4 \i = y_0 e^{\i\phi_4}, \qquad
c_1 + c_3 \i = y_1 e^{\i\phi_4}, \qquad
c_2 + c_6 \i = y_2 e^{\i\phi_4}, \qquad
c_3 + c_7 \i = y_3 e^{\i\phi_4}.
\end{align*}
\medskip

It is worth noticing that the angles $\phi_4,\ldots, \phi_7$ are responsible for the rotations of subsequent (dependent) complex numbers, so that after the last rotation one gets real numbers. Further considerations that could lead to direct formulas for angles $\phi_4,\ldots,\phi_7$ would, however, require a much more complex algebraic apparatus and are currently outside the area of interest of this work.
\medskip

\section{Numerical experiments -- results and problems}
To check what are (experimental) relationships between the angles $\phi_4,\ldots,\phi_7$, the numbers $y_0, \ldots, y_3$, and oktonion $o\in \O$ given by the formula \eqref{eq:oct}, we conducted a series of numerical tests. As in the previous experiment, we generated unitary octonions using spherical coordinates with seven angles $\psi_1,\ldots,\psi_7$, where $\psi_1,\ldots,\psi_6\in[-\pi/2,\pi/2)$ and $\psi_7\in[-\pi,\pi)$. The tests were performed in MATLAB, and the built-in function \texttt{fsolve} was used to solve the given system of equations.
\medskip

\begin{figure}[H]
\centering \includegraphics[width = 0.8\textwidth]{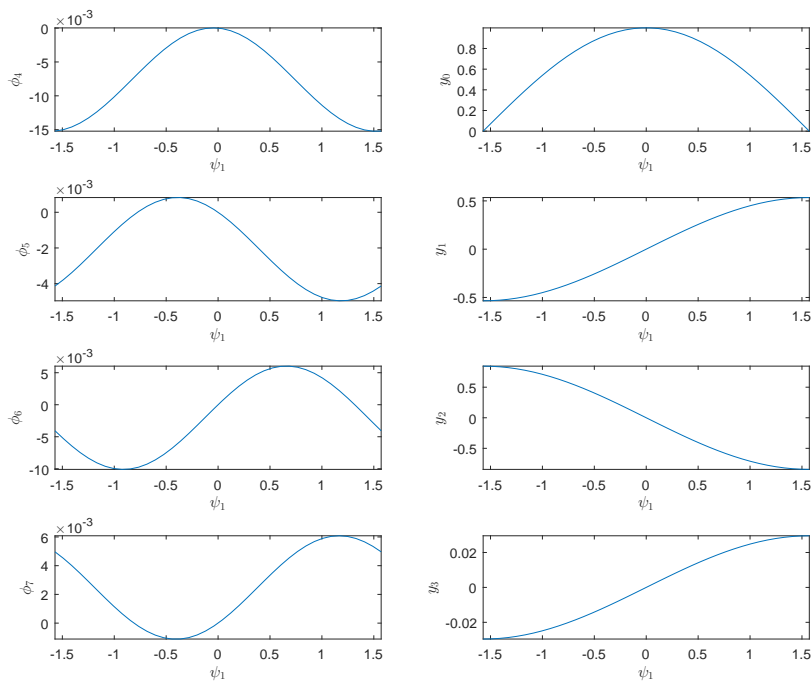}
\caption{Reconstructed angles $\phi_4,\ldots,\phi_7$ and numbers $y_0,\ldots,y_3$ when changing parameter $\psi_1$ and fixed (random) parameters $\psi_2,\ldots,\psi_7$.}
\label{fig:1}
\end{figure}

\begin{figure}[H]
\centering \includegraphics[width = 0.8\textwidth]{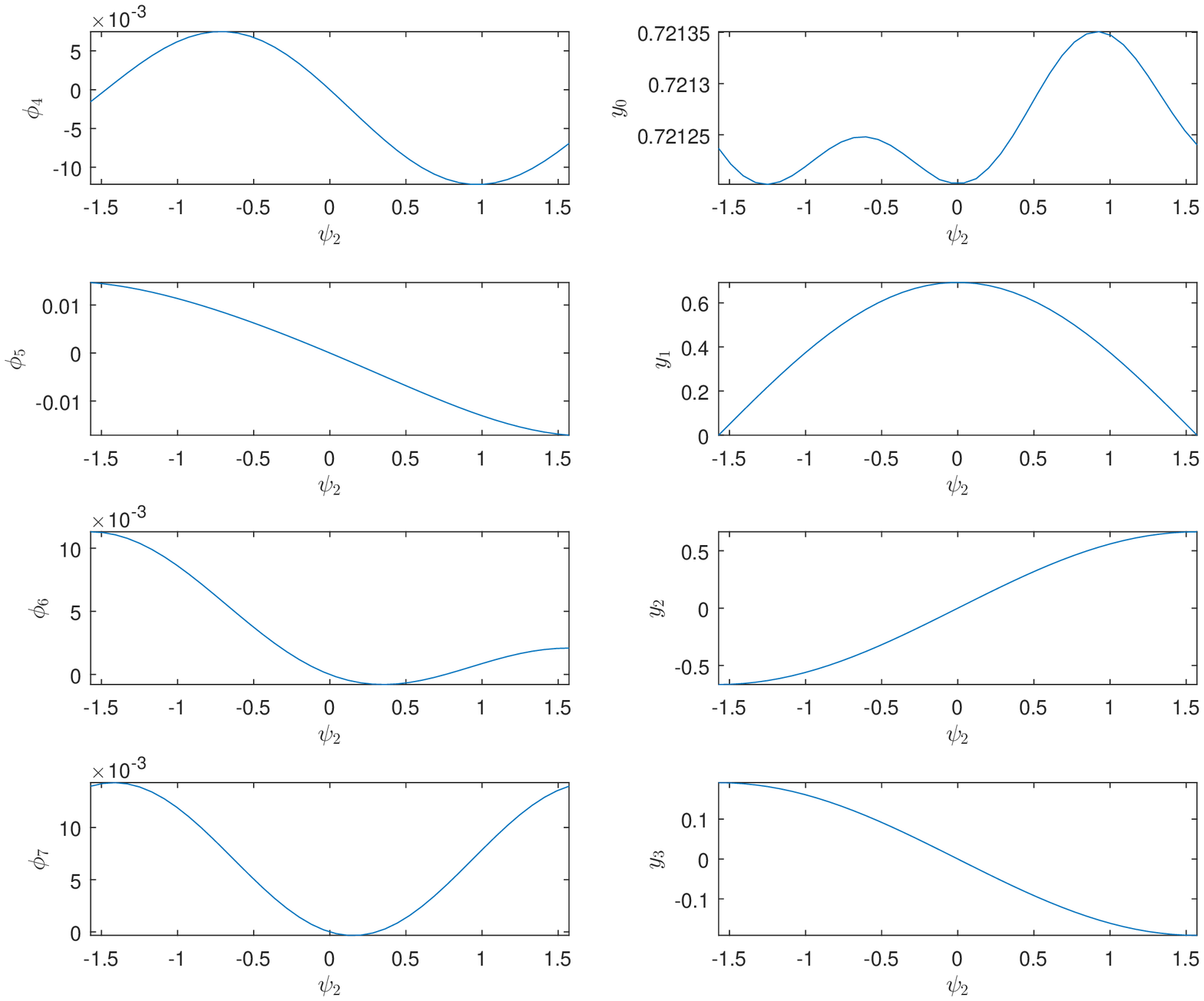}
\caption{Reconstructed angles $\phi_4,\ldots,\phi_7$ and numbers $y_0,\ldots,y_3$ when changing parameter $\psi_2$ and fixed (random) parameters $\psi_1,\psi_3,\ldots,\psi_7$.}
\label{fig:2}
\end{figure}

\begin{figure}[H]
\centering \includegraphics[width = 0.8\textwidth]{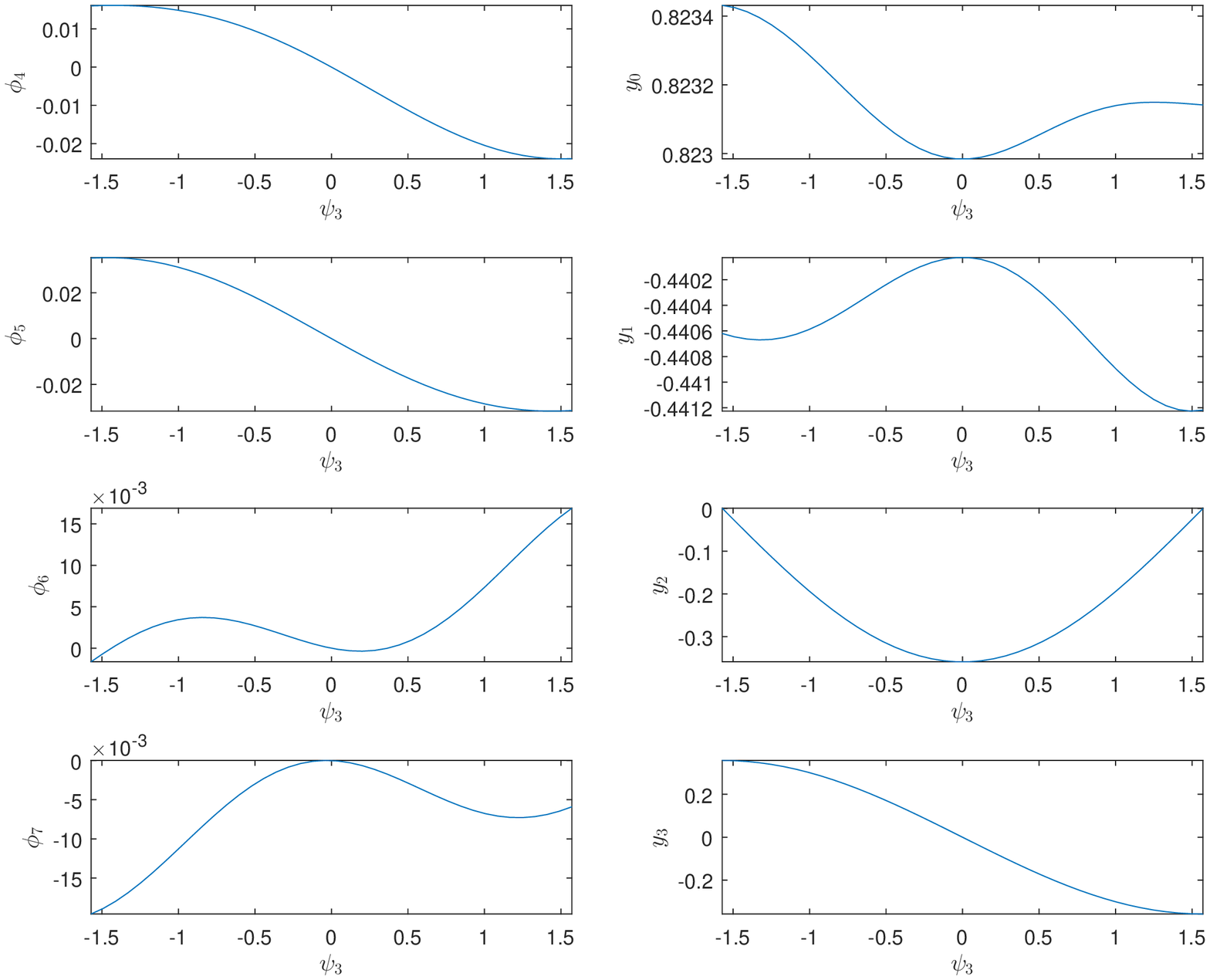}
\caption{Reconstructed angles $\phi_4,\ldots,\phi_7$ and numbers $y_0,\ldots,y_3$ when changing parameter $\psi_3$ and fixed (random) parameters $\psi_1,\psi_2,\psi_4,\ldots,\psi_7$.}
\label{fig:3}
\end{figure}

\begin{figure}[H]
\centering \includegraphics[width = 0.8\textwidth]{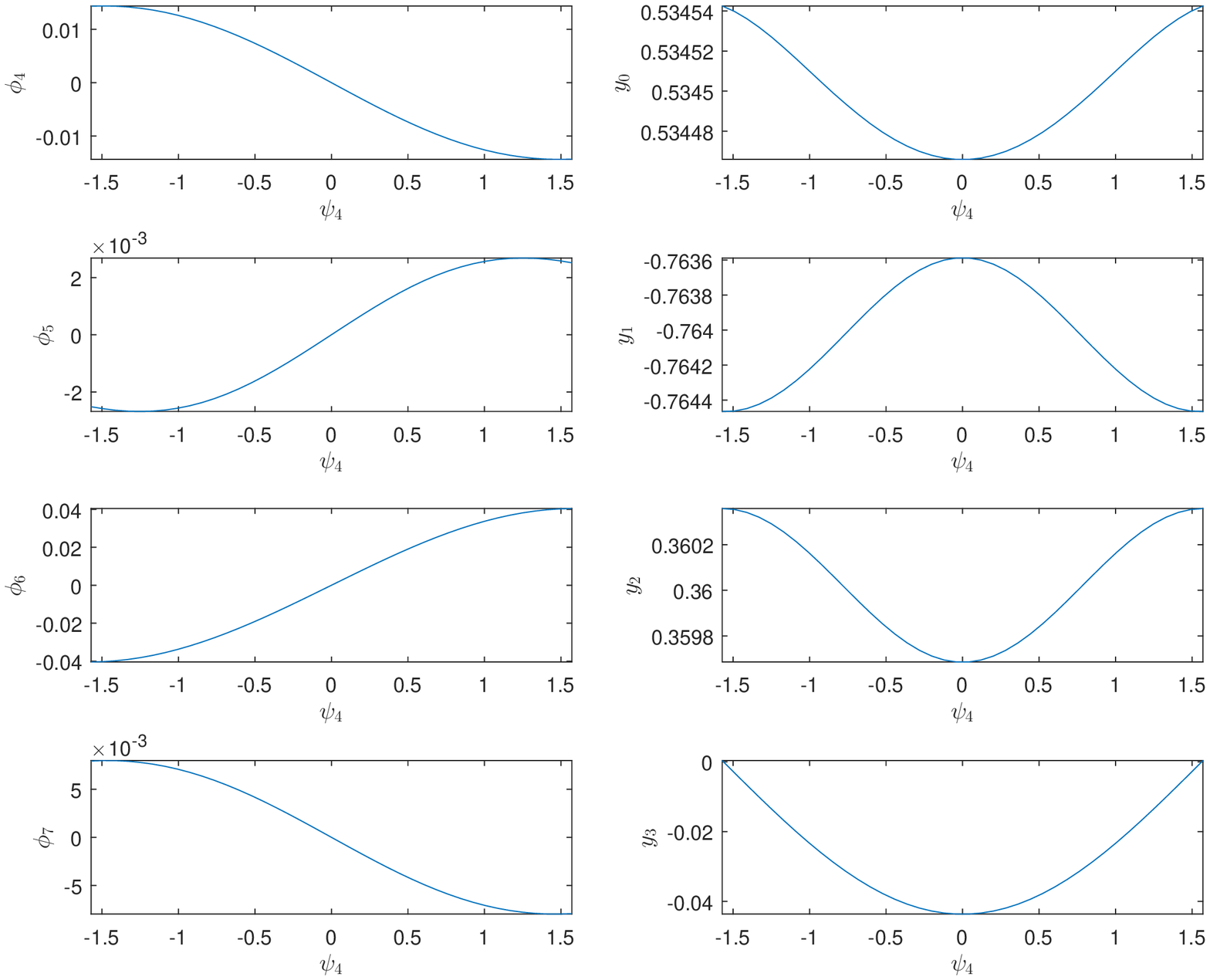}
\caption{Reconstructed angles $\phi_4,\ldots,\phi_7$ and numbers $y_0,\ldots,y_3$ when changing parameter $\psi_4$ and fixed (random) parameters $\psi_1,\ldots,\psi_3,\psi_5,\ldots,\psi_7$.}
\label{fig:4}
\end{figure}

\begin{figure}[H]
\centering \includegraphics[width = 0.8\textwidth]{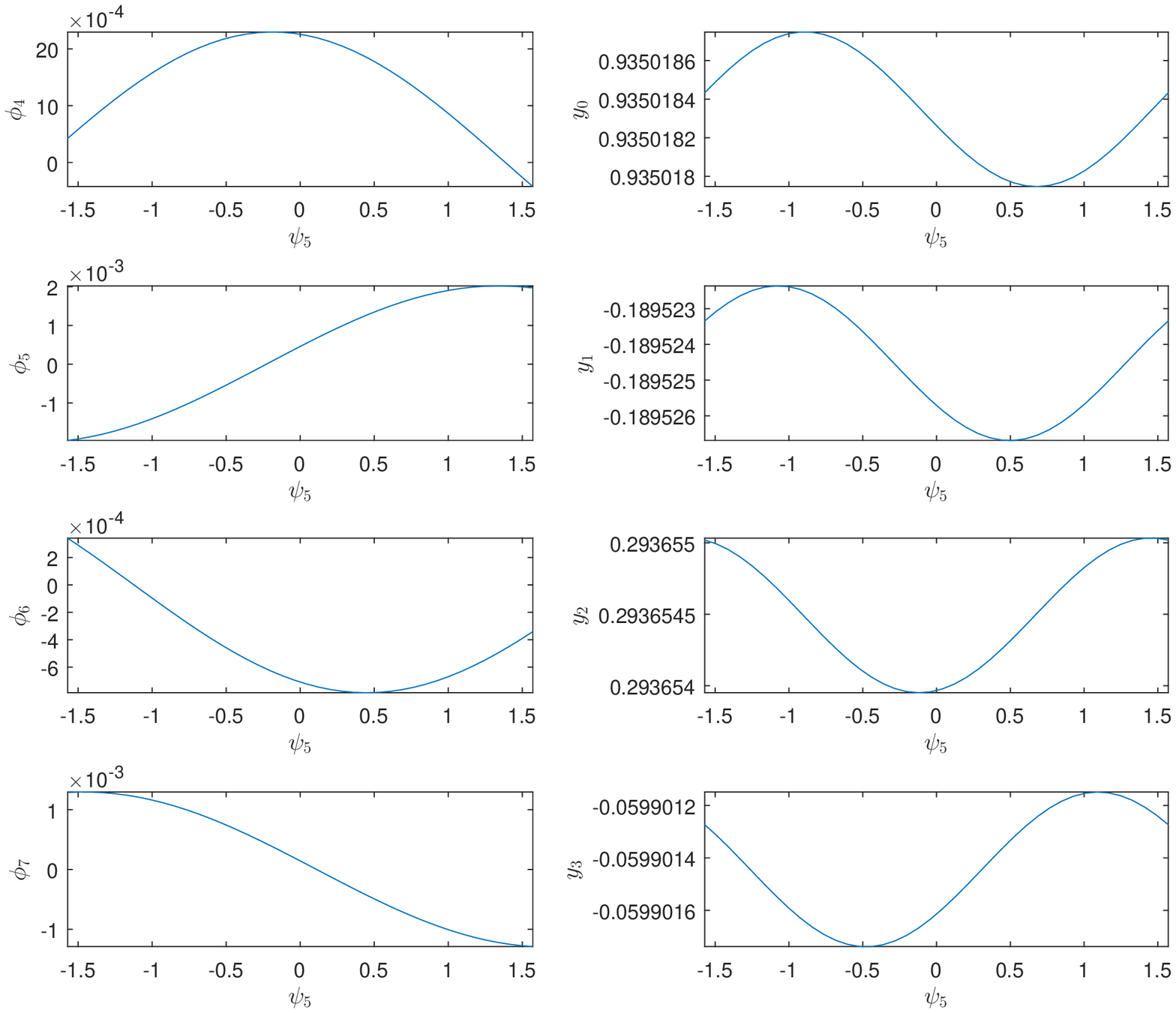}
\caption{Reconstructed angles $\phi_4,\ldots,\phi_7$ and numbers $y_0,\ldots,y_3$ when changing parameter $\psi_5$ and fixed (random) parameters $\psi_1,\ldots,\psi_4,\psi_6,\psi_7$.}
\label{fig:5}
\end{figure}

\begin{figure}[H]
\centering \includegraphics[width = 0.8\textwidth]{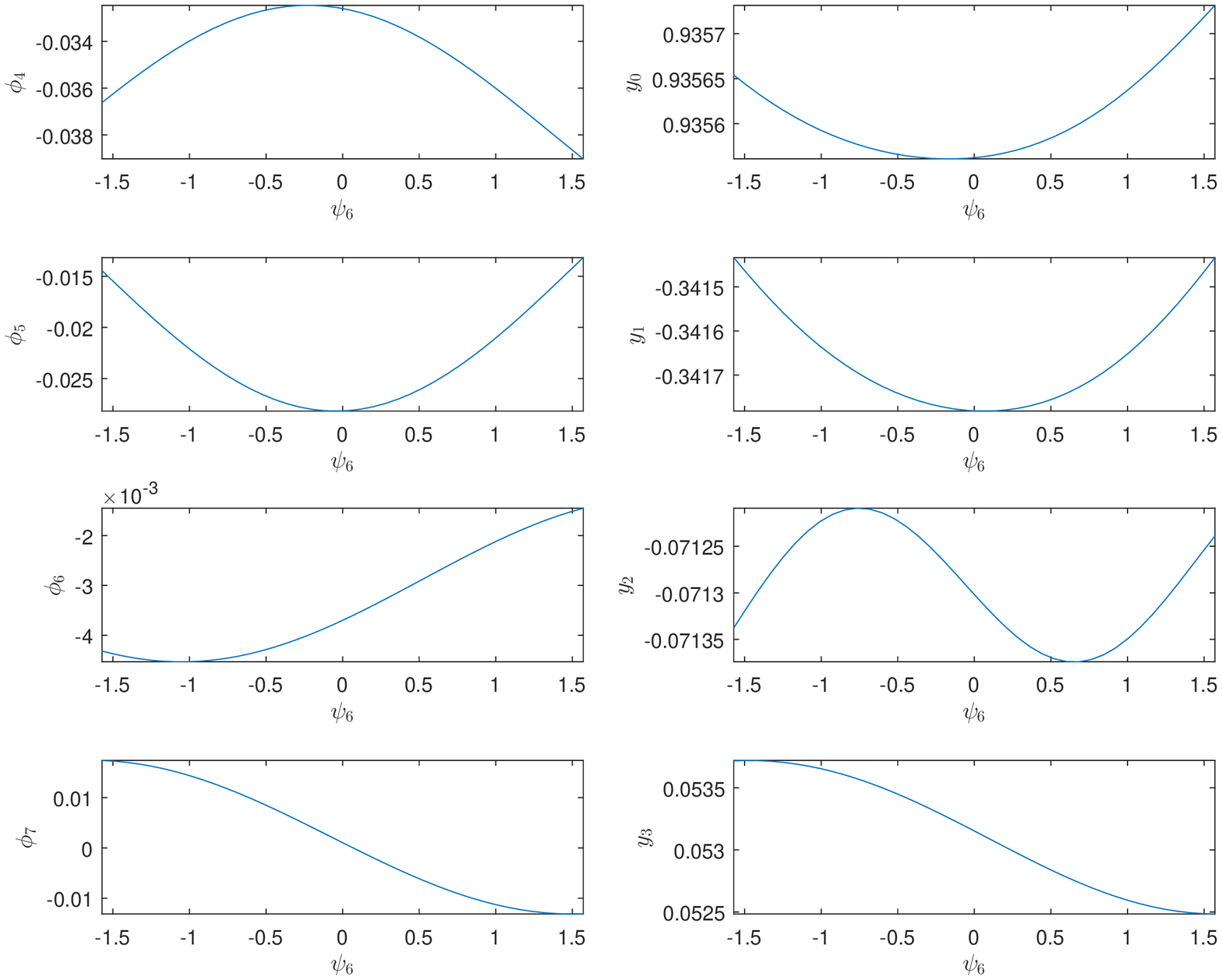}
\caption{Reconstructed angles $\phi_4,\ldots,\phi_7$ and numbers $y_0,\ldots,y_3$ when changing parameter $\psi_6$ and fixed (random) parameters $\psi_1,\ldots,\psi_5,\psi_7$.}
\label{fig:6}
\end{figure}

\begin{figure}[H]
\centering \includegraphics[width = 0.8\textwidth]{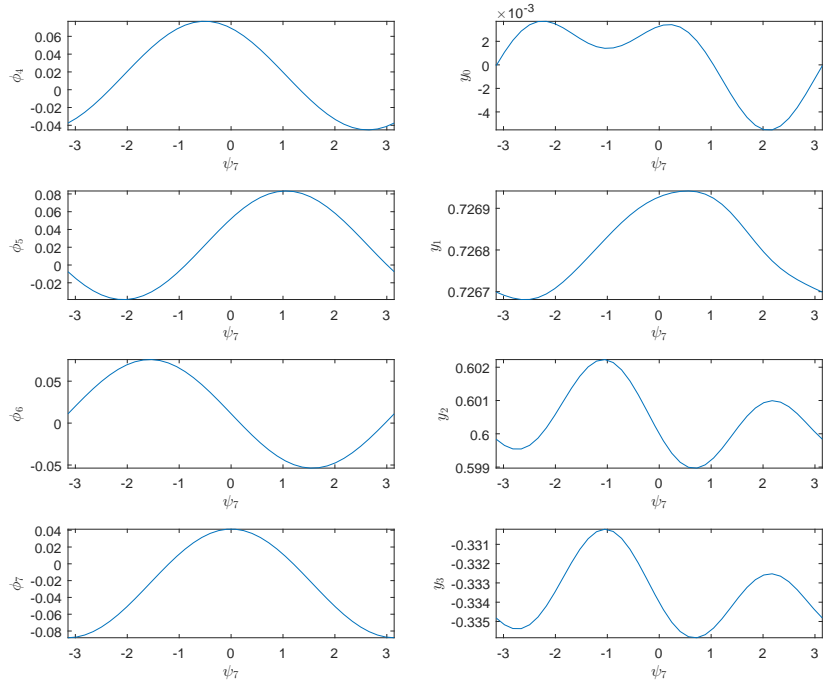}
\caption{Reconstructed angles $\phi_4,\ldots,\phi_7$ and numbers $y_0,\ldots,y_3$ when changing parameter $\psi_7$ and fixed (random) parameters $\psi_1,\ldots,\psi_6$.}
\label{fig:7}
\end{figure}

In Fig. \ref{fig:1}--\ref{fig:7} the results of the experiment are presented. It is worth noting that the obtained dependencies seem to be regular, which suggests that it is possible to find direct formulas describing the relationships between polar representation angles and individual octonion coefficients. Doubts are only aroused by the fact that some of the angles are not periodic functions of the parameters $\psi_1, \ldots, \psi_7$. It can be, however, noticed by extending the ranges of $ \psi_1, \ldots, \psi_6$ to the $ [- \pi, \pi) $ range, that these dependencies then become periodic. It can be deduced from this (and it is not surprising) that the polar representation of octonions is not unique.

\section{Discussion and conclusions}
The issue of the polar (exponential) form of octonion remains almost uncharted at the theoretical angle. The results so far have turned out to be wrong, and the proposed methods only give numerical results. They suggest, however, that analytical formulas can be found, but additional algebraic analysis of the problem is needed.
\medskip

Research on polar representation seems extremely important due to potential applications. In their articles, Hahn and Snopek suggest the possibility of using information about angles in the analysis of hypercomplex analytic signals. On the other hand, knowledge of the angles in the exponential representation would allow a better geometric interpretation of the octonion algebra itself.
\medskip


\subsection*{Acknowledgment}
The~research conducted by the author was supported by National Science Centre (Poland) grant No. 2016/23/N/ST7/00131.

\end{document}